\documentclass[12pt]{article}
\usepackage{amssymb,amsmath,amscd, amsthm,stmaryrd,verbatim, mathrsfs}
\numberwithin{equation}{section}

\begin{document}
\baselineskip=14pt

\newcommand{\la}{\langle}
\newcommand{\ra}{\rangle}
\newcommand{\psp}{\vspace{0.4cm}}
\newcommand{\pse}{\vspace{0.2cm}}
\newcommand{\ptl}{\partial}
\newcommand{\dlt}{\delta}
\newcommand{\sgm}{\sigma}
\newcommand{\al}{\alpha}
\newcommand{\be}{\beta}
\newcommand{\G}{\Gamma}
\newcommand{\gm}{\gamma}
\newcommand{\vs}{\varsigma}
\newcommand{\Lmd}{\Lambda}
\newcommand{\lmd}{\lambda}
\newcommand{\td}{\tilde}
\newcommand{\vf}{\varphi}
\newcommand{\yt}{Y^{\nu}}
\newcommand{\wt}{\mbox{wt}\:}
\newcommand{\rd}{\mbox{Res}}
\newcommand{\ad}{\mbox{ad}}
\newcommand{\stl}{\stackrel}
\newcommand{\ol}{\overline}
\newcommand{\ul}{\underline}
\newcommand{\es}{\epsilon}
\newcommand{\dmd}{\diamond}
\newcommand{\clt}{\clubsuit}
\newcommand{\vt}{\vartheta}
\newcommand{\ves}{\varepsilon}
\newcommand{\dg}{\dagger}
\newcommand{\tr}{\mbox{Tr}}
\newcommand{\ga}{{\cal G}({\cal A})}
\newcommand{\hga}{\hat{\cal G}({\cal A})}
\newcommand{\Edo}{\mbox{End}\:}
\newcommand{\for}{\mbox{for}}
\newcommand{\kn}{\mbox{ker}}
\newcommand{\Dlt}{\Delta}
\newcommand{\rad}{\mbox{Rad}}
\newcommand{\rta}{\rightarrow}
\newcommand{\mbb}{\mathbb}
\newcommand{\lra}{\Longrightarrow}
\newcommand{\X}{{\cal X}}
\newcommand{\Y}{{\cal Y}}
\newcommand{\Z}{{\cal Z}}
\newcommand{\U}{{\cal U}}
\newcommand{\V}{{\cal V}}
\newcommand{\W}{{\cal W}}
\newcommand{\sta}{\theta}
\setlength{\unitlength}{3pt}
\newcommand{\msr}{\mathscr}
\newcommand{\wht}{\widehat}
\newcommand{\mfk}{\mathfrak}

\begin{center}{\large \bf Prime Numbers of the Form $24m+1$\\ \pse and Erd\"{o}s-Straus Conjecture} 
\end{center}

\vspace{0.2cm}

\begin{center}{\large  Yichuan Cao and Xiaoping Xu\footnote{Corresponding Author}}\end{center}

\begin{center}{
State Key Laboratory of Mathematical Sciences, \\ Academy of Mathematics \& System
Sciences\\ Chinese Academy of Sciences, Beijing 100190, P.R. China
\\ \& School of Mathematics, University of Chinese Academy of Sciences,\\ Beijing 100049, P.R. China}\end{center}

\begin {abstract}
\quad

In 1948, Erd\"{o}s and Straus formulated a conjecture on Egyptian fractions. One only needs to prove the conjecture for any prime number $n$ such that $n\equiv 1\;(\mbox{mod}\;24)$.
In this article, we give a necessary and sufficient polynomial condition for a number of the form $24m+1$ to be a prime number. Then we present  56 solvability polynomials consisting of our earlier tame solutions, wild  solutions and the polynomials for composite numbers of the form $24m+1$. Numeric test shows that they imply the conjecture. Moreover, we give a list of 36 exponential-polynomial functions
consisting of our earlier 28 tame solutions and the polynomials for composite numbers, whose numeric test shows that Howerton's conjecture on fourteen wild primes holds. Furthermore, we find a version of  Goldbach conjecture modulo 24.

\psp

\noindent{\it Keywords}:\hspace{0.3cm} Erd\"{o}s-Straus conjecture; Egyptian fraction; prime number; Goldbach conjecture; Howerton conjecture.
\end{abstract}

\section {Introduction}

Ancient Egyptians used sums of unit fractions (whose numerators are 1) to express fractions due to their ways of distributing food.
For instance,  $5/8$ was interpreted by ancient Egyptians as distributing five pancakes among eight people fairly. They cut first four
pancakes into halves and cut the last one into eight equal pieces. So each person got the same share: $1/2+1/8$ pancakes. This amazingly interpreted the mathematical equation
\begin{equation}\frac{5}{8}=\frac{1}{2}+\frac{1}{8}.\end{equation}
So an Egyptian fraction is a sum of distinct unit fractions. By repeatedly applying the simple fact
\begin{equation}\frac{1}{k}=\frac{1}{k+1}+\frac{1}{k(k+1)},\end{equation}
one can easily prove that any fraction is a Egyptian fraction. However,  it is very difficult to determine if a fraction can be expressed as a sum of fixed number of unit fractions (cf. \cite{ETao} for an excellent exposition and extensive references). For example, it is difficult to know if a fraction can be written as a sum of two unit fractions (e.g., cf. \cite{HV1,HV2,HV3,Jc2, Ps}). In 1948, Erd\"{o}s and Straus formulated a conjecture : for any positive integer $n>2$, there exist positive integers $n_1,n_2$ and $n_3$ such that
\begin{equation}\frac{4}{n}=\frac{1}{n_1}+\frac{1}{n_2}+\frac{1}{n_3},\end{equation}
which is still open up to now. For each $n$, the number of solution can go to infinity as $n$ does. Elscholtz and Tao \cite{ETao} found excellent bounds for it.

If $n=24m+1$ and $n_1\leq n_2,n_3$, then $n_1=6m+k$ with $1\leq k\leq 12m$. A solution $(n_1,n_2,n_3)$ of the equation (1.3) is called a {\it tame solution} if $n_2$ and $n_3$ are factors of $(6m+k)(24m+1)$. We call $n=24m+1$ {\it wild} if it does not have any tame solution.
In \cite{Xx1}, we derived
twenty-eight families of tame solutions of the equation (1.3). Numeric test shows
that  they  cover all the tame primes. Based on the information in our earlier
 work on tame solutions posed in arXiv, Howerton \cite{Hb} found that there are only fourteen wild primes of the form $n=24m+1\leq 2.4\times 10^{11}$. In \cite{Xx2}. we derived thirty-four families of wild solutions of the  equation in (1.3), which contain the solvability of the fourteen wild primes. Together with our earlier tame polynomial solutions in \cite{Xx1}, numeric test shows that they cover all the primes of the form $24m+1$.

In this article, we give a necessary and sufficient polynomial condition for a number of the form $24m+1$ to be a prime number. Then we present the 56 solvability polynomials consisting of our earlier tame polynomial solutions in \cite{Xx1}, wild polynomial solutions in \cite{Xx2} and the polynomials in this paper for composite numbers of the form $24m+1$. Numeric test shows that they imply  Erd\"{o}s-Straus conjecture. Moreover, we give a list of 36 exponential-polynomial functions  consisting of our earlier 28 tame solutions in \cite{Xx1} and the polynomials for composite numbers, whose numeric test shows that Howerton's conjecture on fourteen wild primes holds.

Goldbach conjecture says that any even integer $\geq 4$ is a sum of two primes, which is still open today. As a byproduct of this work,
we find a version of  Goldbach conjecture modulo 24. The paper is organized as follows.\pse

In Section 2, we derive a necessary and sufficient condition for a number of the form $24m+1$ to be a prime. In Section 3, we present 56
solvability polynomials for Erd\"{o}s-Straus conjecture. In Section 4, we present 36
solvability  exponential-polynomial functions for Howerton's conjecture. Moreover, we given Goldbach conjecture modulo 24 in Section 5.

\section{Primes of the Form $24m+1$}

In this section, we want to find the conditions that a positive integer of the form $24m+1$ is a prime.

Let $p>3$ be a prime. Then one of the following conditions must hold:
\begin{equation}(1)\quad p\equiv 1\;\;(\mbox{mod}\:12);\quad (2)\quad p\equiv 5\;\;(\mbox{mod}\:12);\end{equation}
\begin{equation}(3)\quad p\equiv 7\;\;(\mbox{mod}\:12);\quad
(4)\quad p\equiv 11\;\;(\mbox{mod}\:12).\end{equation}
In Case (1), $p=12k+1$ for some $k\in\mbb N$, and
\begin{equation}(12k+1)^2=144k^2+24k+1=24k(6k+1)+1.\end{equation}
In Case (2), $p=12k+5$ for some $k\in\mbb N$, and
\begin{equation}(12k+5)^2=144k^2+120k+25=24(6k^2+5k+1)+1=24(2k+1)(3k+1)+1.\end{equation}
In Case (3), $p=12k+7$ for some $k\in\mbb N$, and
\begin{equation}(12k+7)^2=144k^2+168k+49=24(6k^2+7k+2)+1=24(2k+1)(3k+2)+1.\end{equation}
In Case (4), $p=12k+11$ for some $k\in\mbb N$, and
\begin{equation}(12k+11)^2=144k^2+264k+121=24(6k^2+11k+5)+1=24(k+1)(6k+5)+1.\end{equation}

Consider $n=24m+1$ with $m\in\mbb N$. Since $g.c.d(24,p)=1$, there exists $\xi\in\mbb N$ such that
\begin{equation}24\xi\equiv 1\;\;(\mbox{mod}\:p).\end{equation} In particular,
\begin{equation}n\equiv 0\;\;(\mbox{mod}\:p)\lra m\equiv -\xi\;\;(\mbox{mod}\:p).\end{equation}
In Case (1), (2.3) gives $\xi=-k(6k+1)$ and so
\begin{equation}24m+1\equiv 0\;\;(\mbox{mod}\:12k+1)\Leftrightarrow m\equiv k(6k+1)\;\;(\mbox{mod}\:12k+1).\end{equation}
Indeed the above equation holds for any $1<k\in\mbb N$. Similarly, we obtain that
\begin{equation}24m+1\equiv 0\;\;(\mbox{mod}\:12k+5)\Leftrightarrow m\equiv (2k+1)(3k+1)\;\;(\mbox{mod}\:12k+5),\end{equation}
\begin{equation}24m+1\equiv 0\;\;(\mbox{mod}\:12k+7)\Leftrightarrow m\equiv (2k+1)(3k+2)\;\;(\mbox{mod}\:12k+7)\end{equation}and
\begin{equation}24m+1\equiv 0\;\;(\mbox{mod}\:12k+11)\Leftrightarrow m\equiv (k+1)(6k+5)\;\;(\mbox{mod}\:12k+11)\end{equation}
hold for any $k\in\mbb N$.

Next we want to show that the quadratic conditions in the above four expressions can be simplified to linear conditions. \pse

{\it Case (1a)}.\quad Assume $k=2\ell$ in (2.9) with $0<\ell\in\mbb N$.\pse

In this case,
\begin{equation}k(6k+1)=2\ell(12\ell+1)=\ell(24\ell+1)+\ell.\end{equation}Thus
\begin{equation}24m+1\equiv 0\;\;(\mbox{mod}\:24\ell+1)\Leftrightarrow m\equiv \ell\;\;(\mbox{mod}\:24\ell+1).\end{equation}\pse

{\it Case (1b)}.\quad Assume $k=2\ell+1$ in (2.9) with $\ell\in\mbb N$.\pse

In this case,
\begin{equation}k(6k+1)=(2\ell+1)(12\ell+7)=\ell(24\ell+13)+13\ell+7.\end{equation}Thus
\begin{equation}24m+1\equiv 0\;\;(\mbox{mod}\:24\ell+13)\Leftrightarrow m\equiv 13\ell+7\;\;(\mbox{mod}\:24\ell+13).\end{equation}\pse

{\it Case (2a)}.\quad Assume $k=2\ell$ in (2.10) with $\ell\in\mbb N$.\pse

In this case,
\begin{equation}(2k+1)(3k+1)=(4\ell+1)(6\ell+1)=24\ell^2+10\ell+1=\ell(24\ell+5)+5\ell+1.\end{equation}Hence
\begin{equation}24m+1\equiv 0\;\;(\mbox{mod}\:24\ell+5)\Leftrightarrow m\equiv 5\ell+1\;\;(\mbox{mod}\:24\ell+5).\end{equation}\pse

{\it Case (2b)}.\quad Assume $k=2\ell+1$ in (2.10) with $\ell\in\mbb N$.\pse

In this case,
\begin{equation}(2k+1)(3k+1)=(4\ell+3)(6\ell+4)=24\ell^2+34\ell+12=\ell(24\ell+17)+17\ell+12.\end{equation}Thus
\begin{equation}24m+1\equiv 0\;\;(\mbox{mod}\:24\ell+17)\Leftrightarrow m\equiv 17\ell+12\;\;(\mbox{mod}\:24\ell+17).\end{equation}\pse

{\it Case (3a)}.\quad Assume $k=2\ell$ in (2.11) with $\ell\in\mbb N$.\pse

In this case,
\begin{equation} (2k+1)(3k+2)=(4\ell+1)(6\ell+2)=24\ell^2+14\ell+2=\ell(24\ell+7)+7\ell+2.\end{equation}Hence
\begin{equation}24m+1\equiv 0\;\;(\mbox{mod}\:24\ell+7)\Leftrightarrow m\equiv 7\ell+2\;\;(\mbox{mod}\:24\ell+7).\end{equation}\pse

{\it Case (3b)}.\quad Assume $k=2\ell+1$ in (2.11) with $\ell\in\mbb N$.\pse

In this case,
\begin{equation}(2k+1)(3k+2)=(4\ell+3)(6\ell+5)=24\ell^2+38\ell+15=\ell(24\ell+19)+19\ell+15.\end{equation}Thus
\begin{equation}24m+1\equiv 0\;\;(\mbox{mod}\:24\ell+19)\Leftrightarrow m\equiv 19\ell+15\;\;(\mbox{mod}\:24\ell+19).\end{equation}\pse

{\it Case (4a)}.\quad Assume $k=2\ell$ in (2.12) with $\ell\in\mbb N$.\pse

In this case,
\begin{equation} (k+1)(6k+5)=(2\ell+1)(12\ell+5)=24\ell^2+22\ell+5=\ell(24\ell+11)+11\ell+5.\end{equation}Hence
\begin{equation}24m+1\equiv 0\;\;(\mbox{mod}\:24\ell+11)\Leftrightarrow m\equiv 11\ell+5\;\;(\mbox{mod}\:24\ell+11).\end{equation}\pse

{\it Case (4b)}.\quad Assume $k=2\ell+1$ in (2.12) with $\ell\in\mbb N$.\pse

In this case,
\begin{equation}(k+1)(6k+5)=(2\ell+2)(12\ell+11)=24\ell^2+46\ell+22=\ell(24\ell+23)+23\ell+22.\end{equation}Thus
\begin{equation}24m+1\equiv 0\;\;(\mbox{mod}\:24\ell+23)\Leftrightarrow m\equiv 23\ell+22\;\;(\mbox{mod}\:24\ell+23).\end{equation}\pse

In summary, we have:\pse

{\bf Theorem 2.1}\quad {\it Let $m\in\mbb N$. Then $n=24m+1$ is a prime if and only the following condition hold:

(1)\quad $m\not\equiv s\;\;(\mbox{mod}\:24s+1)$ for $0,m\neq s\in\mbb N$;

(2)\quad $m\not\equiv 5s+1\;\;(\mbox{mod}\:24s+5)$ for any $s\in\mbb N$;

(3)\quad $m\not\equiv 7s+2\;\;(\mbox{mod}\:24s+7)$ for any $s\in\mbb N$;

(4)\quad $m\not\equiv 11s+5\;\;(\mbox{mod}\:24s+11)$ for any $s\in\mbb N$;

(5)\quad $m\not\equiv 13s+7\;\;(\mbox{mod}\:24s+13)$ for any $s\in\mbb N$;

(6)\quad $m\not\equiv 17s+12\;\;(\mbox{mod}\:24s+17)$ for any $s\in\mbb N$;

(7)\quad $m\not\equiv 19s+15\;\;(\mbox{mod}\:24s+19)$ for any $s\in\mbb N$;

(8)\quad $m\not\equiv 23s+22\;\;(\mbox{mod}\:24s+23)$ for any $s\in\mbb N$.}\psp

The sufficiency follows from the following proposition.\psp

{\bf Proposition 2.2}. {\it For any positive integer $m$, $n=24m+1$ is a composite number if and only $n=(24s+1)(24t+1)$ for some positive integers $s$ and $t$; or $n=(24s+p)(24t+p)$ for some $p\in\{5,7,11,13,17,19,23\}$ and $s,t\in \mbb N$.}\psp

{\it Proof}. If (1) in the above theorem fails,
 then
\begin{equation}m=s+t(24s+1)=s+t+24st\qquad\for\;\;s,t\in\mbb N,\;st\neq 9.\end{equation}
When (2) in the above theorem fails, we have
\begin{equation}m=5s+1+t(24s+5)=1+5(s+t)+24st\qquad\for\;\;s,t\in\mbb N,\;st\neq 9.\end{equation}
In the case that (3) in the above theorem fails, we have
\begin{equation}m=7s+2+t(24s+7)=2+7(s+t)+24st\qquad\for\;\;s,t\in\mbb N,\;st\neq 9.\end{equation}
Suppose that (4) in the above theorem fails. We have
\begin{equation}m=11s+5+t(24s+11)=5+11(s+t)+24st\qquad\for\;\;s,t\in\mbb N,\;st\neq 9.\end{equation}
If (5) in the above theorem fails, then
\begin{equation}m=13s+7+t(24s+13)=7+13(s+t)+24st\qquad\for\;\;s,t\in\mbb N,\;st\neq 9.\end{equation}
When (6) in the above theorem fails, we have
\begin{equation}m=17s+12+t(24s+17)=12+17(s+t)+24st\qquad\for\;\;s,t\in\mbb N,\;st\neq 9.\end{equation}
In the case that (7) in the above theorem fails, we have
\begin{equation}m=19s+15+t(24s+19)=15+19(s+t)+24st\qquad\for\;\;s,t\in\mbb N,\;st\neq 9.\end{equation}
Suppose that (8) in the above theorem fails. We have
\begin{equation}m=23s+22+t(24s+23)=22+23(s+t)+24st\qquad\for\;\;s,t\in\mbb N,\;st\neq 9.\end{equation}

Note
\begin{equation}24+1=5^2,\quad 24\times 2+1=7^2,\quad  24\times 5+1=121=11^2,\quad 24\times 7+1=169=13^2, \end{equation}
\begin{equation} 24\times 12+1=289=17^2,\quad 24\times 15+1=361=19^2, \quad 24\times 22+1=529=23^2.\end{equation}
Thus we have
\begin{equation}24(s+t+24st)+1=(24s+1)(24t+1),\end{equation}
\begin{equation}24[1+5(s+t)+24st]+1=24\times 1+1+24[5(s+t)+24st]=(24s+5)(24t+5),\end{equation}
\begin{equation}24[2+7(s+t)+24st]+1=24\times 2+1+24[7(s+t)+24st]=(24s+7)(24t+7),\end{equation}
\begin{equation}24[5+11(s+t)+24st]+1=24\times 5+1+24[11(s+t)+24st]=(24s+11)(24t+11),\end{equation}
\begin{equation}24[7+13(s+t)+24st]+1=24\times 7+1+24[13(s+t)+24st]=(24s+13)(24t+13),\end{equation}
\begin{equation}24[12+17(s+t)+24st]+1=24\times 12+1+24[17(s+t)+24st]=(24s+17)(24t+17),\end{equation}
\begin{equation}24[15+19(s+t)+24st]+1=24\times 15+1+24[19(s+t)+24st]=(24s+19)(24t+19),\end{equation}
\begin{equation}24[22+23(s+t)+24st]+1=24\times 22+1+24[23(s+t)+24st]=(24s+23)(24t+23).\end{equation}\pse

Now the conclusion follows from (2.29)-(2.36) and (2.39)-(2.46). $\qquad\Box$\psp

We set\pse

\begin{center}{\bf \large Table 2.1}\quad Polynomials for Composite Numbers \end{center}

\begin{center}\begin{tabular}{|c||c||c||c|}\hline $i$&$f_i(x,y)$&$i$&$f_i(x,y)$\\\hline
1&$x+y+24xy$&2&$1+5(x+y)+24xy$\\\hline
3&$2+7(x+y)+24xy$&4&$5+11(x+y)+24xy$\\\hline
5&$7+13(x+y)+24xy$&6&$12+17(x+y)+24xy$\\\hline
7&$15+19(x+y)+24xy$&8&$22+23(x+y)+24xy$\\\hline
\end{tabular}\end{center}\psp

 Denote by $\mbb Z_+$ the set of positive integers. \psp

 {\bf Corollary 2.3}. {\it A positive integer of the form
 $24m+1$ is composite if and only if $m=f_1(x,y)$ for $x,y\in\mbb Z_+$, or $m=f_i(x,y)$ for some $i\in\{2,3,4,5,6,7,8\}$ and $x,y\in\mbb N$.}  \psp

Set
\begin{equation}{\msr S}_1=\{f_1(x,y)\mid x,y\in\mbb Z_+\},\end{equation}
\begin{equation}{\msr S}_i=\{f_i(x,y)\mid x,y\in\mbb N\}\end{equation}
for $2\leq i\leq 8$. Let
\begin{equation}{\msr S}=\bigcup_{i=1}^8{\msr S}_i.\end{equation}
\pse

{\bf Corollary 2.4}. {\it A positive integer of the form
 $24m+1$ is a prime if and only $m\not\in {\msr S}$.}\psp

In particular, we have the following table:

\begin{center}{\bf \large Table 2.2}\quad $1\leq m\leq 1000$ for Prime $24m+1$\end{center}

\begin{center}\begin{tabular}{|c|c|c|c|c|c|c|c|c|c|c|c|c|c|c|c|}\hline
3&4&8&10&13&14&17&18&19&24&25&28&32&39&42\\\hline 43& 47&48&50&52&54&55&62& 67&69&73&74&75&78&83\\\hline 84&87& 88&89&90& 95&99&103&105&108&109&112&113&118&119\\\hline 123&125&127& 130&132&134&138&140&143&144&147&153&154&157&158\\\hline 162&169&172&173&174&175&178&179&185&188&190&194&197&200&207\\\hline 208&213&217&218&220&227&230&232&235&237&239&244&245&248&253\\\hline 255&259&264&265&270&272&273&274&278&283&285&290&294&297&299\\\hline 304&305&307&308&309&312&314&315&320&323&328&333&334&337&340\\\hline 342&343&347&348&349&355&360&362&363&364&365&372&375&377&385\\\hline 389&393&400&402&404&405&407&409&417&424&428&430&432&438&444\\\hline
347&348&358&463&465&469&472&473&479&483&484&467&493&498&502\\\hline 503&504&510&512&517&518&519&523&524&525&529&530&535&537&542\\\hline
543&549&552&554&559&560&563&564&568&570&572&578&580&595&600\\\hline 602&608&613&614&622&628&630&633&634&635&637&638&640&650&652\\\hline
659&662&663&664&668&669&677&678&682&684&690&693&694&697&705\\\hline 708&710&714&717&719&724&725&727&729&732&738&739&740&745&747\\\hline
749&752&754&755&757&759&762&763&768&769&770&773&783&788&792\\\hline 795&802&803&809&810&812&817&820&823&824&825&833&837&838&840\\\hline
843&848&855&858&860&867&869&872&875&880&882&883&888&893&895\\\hline 897&899&900&902&903&909&910&914&915&922&923&928&932&935&949\\\hline
955&958&959&960&967&978&979&983&987&990&993&994&997&999&1000\\\hline
\end{tabular}\end{center}

For instance, we have the following primes:
\begin{equation}24\times 3+1=73,\quad 24\times 4+1=96,\quad 24\times 8+1=193,\quad 24\times 10+1=241,\end{equation}
\begin{equation}24\times 13+1=313,\quad 24\times 14+1=337,\quad 24\times 17+1=409,\quad 24\times 18+1=433,\end{equation}
\begin{equation}24\times 19+1=457,\quad 24\times 24+1=577,\quad 24\times 25+1=601,\quad 24\times 28+1=673,\end{equation}
\begin{equation}24\times 32+1=769,\quad 24\times 39+1=937,\quad 24\times 42+1=1009,\quad 24\times 43+1=1033.\end{equation}

\section{Solvability Polynomials for E-S Conjecture}

 The following table is taken from the wild solutions of the Erd\"{o}s-Straus equation (1.3) found in \cite{Xx2}:

\begin{center}{\bf \large Table 3.1\quad Polynomials for Wild Solutions}\end{center}
\begin{center}\begin{tabular}{|c||c||c||c|}\hline $i$&$f_i(x,y)$&$i$&$f_i(x,y)$\\\hline
9&$4+5(x+y)+6xy$&10&$4+10x+7y+15xy$\\\hline
11&$3+7x+10y+21xy$&12&$3+11x+7y+22xy$\\\hline
13&$11+14x+17y+21xy$&14&$18+19y+23x(y+1)$\\\hline
15&$8+20x+13y+30xy$&16&$9+10x+29y+30xy$\\\hline
17&$4+14x+13y+42xy$&18&$3+28x+5y+42xy$\\\hline
19&$6+14x+19y+42xy$&20&$27+28x+41y+42xy$\\\hline
21&$6+7y+46x(y+1)$&22&$14+15y+46x(y+1)$\\\hline
23&$9+33x+19y+66xy$&24&$20+33x+41y+66xy$\\\hline
25&$21+33x+43y+66xy$&26&$32+33x+65y+66xy$\\\hline
27&$9+10y+69x(y+1)$&28&$10+11y+69x(y+1)$\\\hline
29&$33+34y+69x(y+1)$&30&$55+56y+69x(y+1)$\\\hline
31&$73+74y+92x(y+1)$&32&$18+19x+113y+114xy$\\\hline
33&$30+95x+37y+114xy$&34&$16+17y+138x(y+1)$\\\hline
35&$42+43y+138x(y+1)$&36&$66+67y+138x(y+1)$\\\hline
37&$88+89y+138x(y+1)$&38&$90+91y+138x(y+1)$\\\hline
39&$108+109y+138x(y+1)$&40&$112+113y+138x(y+1)$\\\hline
41&$136+137y+138x(y+1)$&&\\\hline

\end{tabular}\end{center}\psp

The details of the above polynomials in \cite{Xx2} are as follows: \\ $m=f_9(x,y)$ is taken from (3.6) and the solution for $n=24m+1$ was given in (3.8); $m=f_{10}(x,y)$ is taken from (3.89) and the solution for $n=24m+1$ was given in (3.91);\\ $m=f_{11}(x,y)$ is taken from (3.26) and the solution for $n=24m+1$ was given in (3.28); $m=f_{12}(x,y)$ is taken from (3.52) and the solution for $n=24m+1$ was given in (3.54);
$m=f_{13}(x,y)$ is taken from (3.32) and the solution for $n=24m+1$ was given in (3.34);\\ $m=f_{14}(x,y)$ is taken from (3.37) and the solution for $n=24m+1$ was given in (3.39);\\ $m=f_{15}(x,y)$ is taken from (3.95) and the solution for $n=24m+1$ was given in (3.97);\\ $m=f_{16}(x,y)$ is taken from (3.83) and the solution for $n=24m+1$ was given in (3.85);\\ $m=f_{17}(x,y)$ is taken from (3.17) and the solution for $n=24m+1$ was given in (3.19);\\ $m=f_{18}(x,y)$ is taken from (3.40) and the solution for $n=24m+1$ was given in (3.42);\\ $m=f_{19}(x,y)$ is taken from (3.45) and the solution for $n=24m+1$ was given in (3.47);\\ $m=f_{20}(x,y)$ is taken from (3.14) and the solution for $n=24m+1$ was given in (3.16);\\
$m=f_{21}(x,y)$ is taken from (3.56) and the solution for $n=24m+1$ was given in (3.58);\\ $m=f_{22}(x,y)$ is taken from (4.5) and the solution for $n=24m+1$ was given in (4.7);\\ $m=f_{23}(x,y)$ is taken from (3.71) and the solution for $n=24m+1$ was given in (3.73);\\ $m=f_{24}(x,y)$ is taken from (3.77) and the solution for $n=24m+1$ was given in (3.79);\\ $m=f_{25}(x,y)$ is taken from (3.54) and the solution for $n=24m+1$ was given in (3.66);\\ $m=f_{26}(x,y)$ is taken from (3.58) and the solution for $n=24m+1$ was given in (3.60);\\ $m=f_{27}(x,y)$ is taken from (3.75) and the solution for $n=24m+1$ was given in (3.77);\\ $m=f_{28}(x,y)$ is taken from (4.43) and the solution for $n=24m+1$ was given in (4.45);\\ $m=f_{29}(x,y)$ is taken from (4.49) and the solution for $n=24m+1$ was given in (4.51);\\ $m=f_{30}(x,y)$ is taken from (4.81) and the solution for $n=24m+1$ was given in (4.83);\\ $m=f_{31}(x,y)$ is taken from (4.88) and the solution for $n=24m+1$ was given in (4.89);\\ $m=f_{32}(x,y)$ is taken from (3.103) and the solution for $n=24m+1$ was given in (3.105);\\ $m=f_{33}(x,y)$ is taken from (3.109) and the solution for $n=24m+1$ was given in (3.111);\\ $m=f_{34}(x,y)$ is taken from (4.98) and the solution for $n=24m+1$ was given in (4.100);\\ $m=f_{35}(x,y)$ is taken from (4.24) and the solution for $n=24m+1$ was given in (4.26);\\ $m=f_{36}(x,y)$ is taken from (4.62) and the solution for $n=24m+1$ was given in (4.64);\\ $m=f_{37}(x,y)$ is taken from (4.29) and the solution for $n=24m+1$ was given in (4.32);\\ $m=f_{38}(x,y)$ is taken from (4.17) and the solution for $n=24m+1$ was given in (4.19);\\ $m=f_{39}(x,y)$ is taken from (4.104) and the solution for $n=24m+1$ was given in (4.106);\\ $m=f_{40}(x,y)$ is taken from (4.68) and the solution for $n=24m+1$ was given in (4.70);\\ $m=f_{41}(x,y)$ is taken from (4.11) and the solution for $n=24m+1$ was given in (4.13). \psp

The following table is taken from the tame solutions of  Erd\"{o}s-Straus equation (1.3) found in \cite{Xx1}:
\begin{center}{\bf \large Table 3.2\quad Polynomials for Tame Solutions}\end{center}

\begin{center}\begin{tabular}{|c||c|}\hline i&$f_i(x,y,z)$\\\hline
42&$x-y+(4x+1)(16y+3+z(24y+5))$\\\hline
43&$(z+1)(4x+3)(24y+5)-x-y-1$\\\hline
44&$7x-y+(16y+8+z(24y+13))(12x+1)$\\\hline
45&$(z+1)(12x+1)(8y+7)-x-y-1$\\\hline
46&$(z+1)(4x+1)(24y+17)-x-y-1$\\\hline
47&$(z+1)(12x+5)(24y+13)-x-y-1$\\\hline
48&$7x-y+4+(12x+7)(16y+z(24y+1))$\\\hline
49&$(z+1)(12x+7)(8y+3)-x-y-1$\\\hline
50&$3x-y+1+(12x+7)(16y+11+z(24y+17))$\\\hline
51&$(z+1)(12x+11)(24y+1)-x-y-1$\\\hline
52&$2x-y+2(16y+7+z(24y+11))(4x+1)$\\\hline
53&$2[(z+1)(4x+3)(24y+11)-x-1]-y$\\\hline
54&$2x-y+2(16y+7+z(24y+11))(4x+3)+1$\\\hline
55&$4(6x+1)(16y+4+z(24y+7))+14x-y+2$\\\hline
56& $5y-2x+2(8x+3+z(12x+5))(8y+3)+1$\\\hline
\end{tabular}\end{center}\psp

The details of the above polynomials in \cite{Xx1} are as follows: \\
$m=f_{42}(x,y,z)$ is taken from (3.15) and the solution for $n=24m+1$ was given in (3.17);\\
$m=f_{43}(x,y,z)$ is taken from (3.24) and the solution for $n=24m+1$ was given in (3.27);\\
$m=f_{44}(x,y,z)$ is taken from (3.38) and the solution for $n=24m+1$ was given in (3.40);\\
$m=f_{45}(x,y,z)$ is taken from (3.47) and the solution for $n=24m+1$ was given in (3.49);\\
$m=f_{46}(x,y,z)$ is taken from (3.56) and the solution for $n=24m+1$ was given in (3.58);\\
$m=f_{47}(x,y,z)$ is taken from (3.65) and the solution for $n=24m+1$ was given in (3.67);\\
$m=f_{48}(x,y,z)$ is taken from (3.78) and the solution for $n=24m+1$ was given in (3.80);\\
$m=f_{45}(x,y,z)$ is taken from (3.90) and the solution for $n=24m+1$ was given in (3.92);\\
$m=f_{50}(x,y,z)$ is taken from (3.103) and the solution for $n=24m+1$ was given in (3.105);\\
$m=f_{51}(x,y,z)$ is taken from (3.114) and the solution for $n=24m+1$ was given in (3.116);\\
$m=f_{52}(x,y,z)$ is taken from (4.16) with $x=s$ and $y=2t$, and the solution for $n=24m+1$ was given in (4.18);\\
$m=f_{53}(x,y,z)$ is taken from (4.24) with $x=s_1$ and $y=2t$, and the solution for $n=24m+1$ was given in (4.26);\\
$m=f_{54}(x,y,z)$ is taken from (4.96) with $x=s$ and $y=2t$, and the solution for $n=24m+1$ was given in (4.98);\\
$m=f_{55}(x,y,z)$ is taken from (5.75) with $r=0,\;x=s$ and $y=4t$, and the solution for $n=24m+1$ was given in (5.77);\\
$m=f_{56}(x,y,z)$ is taken from (4.39) with $x=s$ and $y=2t$, and the solution for $n=24m+1$ was given in (4.41).\pse

Recall that $\mbb N$ is the set of nonnegative integers and $\mbb Z_+$ is the set of positive integers. Take Tables 2.1, 3.1 and 3.2.
Set
\begin{equation}{\msr S}_1=\{f_1(x,y)\mid x,y\in\mbb Z_+\},\end{equation}
\begin{equation}{\msr S}_i=\{f_i(x,y)\mid x,y\in\mbb N\}\end{equation}
for $2\leq i\leq 41$, and
\begin{equation}{\msr S}_j=\{f_j(x,y,z)\mid x,y,z\in\mbb N\}\end{equation}
for $42\leq j\leq 56$.\psp

{\bf Claim}.\quad $\mbb Z_+=\bigcup_{r=1}^{56}{\msr S}_r$.
\psp

We have used AI to test the claim for $m\leq 10^{11}$ and it holds.
But we failed to get  a logic proof. It is known that one only need to prove  Erd\"{o}s-Straus conjecture for a prime of the form $24m+1$ (e.g., cf. Section 2 in \cite{Xx1}). In fact, any solution of the Erd\"{o}s-Straus equation (1.3) for any prime factor in a composite number $n>2$ induces a solution for $n$, Note that $24m+1$ is a primes if and only if $m\not\in \bigcup_{i=1}^8{\msr S}_i$.
The above claim implies
\begin{equation}\{24m+1\mid m\in\mbb Z_+\}=\{24m+1\mid m\in\bigcup_{r=1}^{56}{\msr S}_r\};\end{equation}
equivalently,  Erd\"{o}s-Straus conjecture holds.

\section{Howerton's Conjecture}

In our first version of \cite{Xx1}, we gave the definitions of tame prime and wild prime. Moreover, we found that there are only nine wild primes among $24m+1$ with $m\leq 30000$. Based the above information, Howerton \cite{Hb} found that there are only fourteen wild primes of the form $24m+1$ with $m\leq 10^{10}$. They are
\begin{equation}409,\quad 577,\quad 5569,\quad 9601,\quad 23929,\quad 83449,\quad 102001,\quad 329617,\quad 712321,\end{equation}
 \begin{equation}1134241,\quad 1724209,\quad 1726201,\quad 5212561,\quad 8813281.\end{equation}
He made a conjecture:
\begin{equation}\mbox{\bf There Are Exactly Fourteen Wild Primes !}\end{equation}
In the following, we want to give a theoretic evidence of supporting Howerton's conjecture. Moreover, we redefine the following $f_i$'s:

\begin{center}{\bf \large Table 4.1\quad More Polynomials for Tame Solutions}\end{center}

\begin{center}\begin{tabular}{|c||c|}\hline i&$f_i(x,y,z)$
\\\hline 9&$x-y+(4x+1)(16y+3+z(24y+5))$\\\hline 10&$(z+1)(4x+3)(24y+5)-x-y-1$\\\hline
11&$7x-y+(16y+8+z(24y+13))(12x+1)$\\\hline
12&$(z+1)(12x+1)(8y+7)-x-y-1$\\\hline
13&$(z+1)(4x+1)(24y+17)-x-y-1$\\\hline
14&$(z+1)(12x+5)(24y+13)-x-y-1$\\\hline
15&$7x-y+4+(12x+7)(16y+z(24y+1))$\\\hline
16&$(z+1)(12x+7)(8y+3)-x-y-1$\\\hline
17&$3x-y+1+(12x+7)(16y+11+z(24y+17))$\\\hline
18&$(z+1)(12x+11)(24y+1)-x-y-1$\\\hline
19&$2x-y+2(16y+7+z(24y+11))(4x+1)$\\\hline
20&$2[(z+1)(4x+3)(24y+11)-x-1]-y$\\\hline 21& $5y-2x+2(8x+3+z(12x+5))(8y+3)+1$\\\hline 22&
$2[(z+1)(12x+5)(8y+1)-x]-y-1$\\\hline 23&$2[(z+1)(12x+7)(24y+19)-x-1]-y$\\\hline
24&$2[(z+1)(12x+1)(24y+19)-x]-y-1$\\\hline
25&$2[(z+1)(4x+1)(24y+11)-x]-y-1$\\\hline 26&$2x-y+2(16y+7+z(24y+11))(4x+3)+1$\\\hline
27&$5y+2[(8x+7+z(12x+11))(8y+3)-x]$\\\hline 28& $2[(z+1)(12x+11)(8y+1)-x-1]-y$\\\hline

\end{tabular}\end{center}and

\begin{center}{\bf \large Table 4.2\quad Exponential-Polynomials for Tame Solutions}\end{center}
\begin{center}\begin{tabular}{|c||c|}\hline i &$f_i(r,x,y,z)$\\\hline
29&$2^r(2x+1)[4(16y+15+z(24y+23))+1]-y-1$\\\hline
30&$2^{2r+2}[2(6x+5)(16y+4+z(24y+7))+7x+6]-y-(2^{2r+1}+1)/3$\\\hline
31&$2^{2r+1}[2(6x+1)(16y+4+z(24y+7))+7x+1]-y+(2^{2r}-1)/3$\\\hline
32&$(z+1)2^{r+2}(2x+1)(24y+23)-y-2^{r}(2x+1)-1$\\\hline
33& $(z+1)2^{2r+3}(6x+1)(24y+7)-y-2^{2r+2}x-(2^{2r+1}+1)/3$\\\hline
34&$(z+1)2^{2r+3}(6x+5)(8y+5)-y-2^{2r+2}x-(5\cdot 2^{2r+1}+2)/3$\\\hline
35&$(z+1)2^{2r+2}(6x+1)(8y+5)-y-2^{2r+1}x-(2^{2r}+2)/3$\\\hline
36&$(z+1)2^{2r+2}(6x+5)(24y+7)-y-2^{2r+1}x-(5\cdot 2^{2r}+1)/3$\\\hline
\end{tabular}\end{center}\psp

The details of the above functions in \cite{Xx1} are as follows: \\
the polynomials $f_9(x,y,z)$-$f_{18}(x,y,z)$ in Table 4.1 are exactly $f_{42}(x,y,z)$-$f_{51}(x,y,z)$ in Table 3,2, respectively;
\\ the polynomial $f_{19}(x,y,z)$-$f_{28}(x,y,z)$ are taken from the table in Section 4.7 of \cite{Xx1};\\
the exponential-polynomial functions $f_{29}(r,x,y,z)$-$f_{36}(r,x,y,z)$ are taken from the table in Section 5.9 of \cite{Xx1}.\\

Denote
\begin{eqnarray}{\msr W}&=&\{17,\;24,\;232,\;400,\;997,\;3477,\;4250,\;13734,\;29680,\nonumber \\&& 47260,\;71842,\;71925,\;217190,\;367220\}.\end{eqnarray}
Then $\{24m+1\mid m\in \msr W\}$ are exactly the wild primes in (4.1) and (4.2).

Recall that $\mbb N$ is the set of nonnegative integers and $\mbb Z_+$ is the set of positive integers.
Take Table 2.1 and denote
\begin{equation}{\msr T}_1=\{f_1(x,y)\mid x,y\in\mbb Z_+\};\end{equation}
\begin{equation}{\msr T}_i=\{f_i(x,y)\mid x,y\in\mbb N\}\end{equation}
for $2\leq i\leq 8$. Take Table 4.1 and set
\begin{equation}{\msr T}_i=\{f_i(x,y,z)\mid x,y,z\in\mbb N\}\end{equation}
for $9\leq i\leq 28$. Take Table 4.2 and let
\begin{equation}{\msr T}_i=\{f_i(r,x,y,z)\mid r,x,y,z\in\mbb N\}\end{equation}
for $29\leq i\leq 36$. Define
\begin{equation}{\msr T}=\bigcup_{i=1}^{36}{\msr T}_i. \end{equation}\pse

{\bf Claim}.\quad {\it For any positive integer $m$, if $m\not\in {\msr T}$, then $m\in {\msr W}$}.\psp

We have used AI to test the claim for $m\leq 10^7$ and it holds. But we failed to get a logic proof. In terms of the form $24m+1$,
the above claim implies Howerton's conjecture in (4.3). Moreover, Howerton's conjecture is stronger than  Erd\"{o}s-Straus  conjecture.

\section{Goldbach Conjecture Modolo 24}

In his letter to Euler 1874, Goldbach proposed the following conjecture: any even integer $\geq 4$ is a sum of two primes.
The most advanced progress in this direction is Chen's Theorem which was announced in 1966 and the detailed proof was published in 1973 (cf. \cite{C1,C2,C3}). Chen's Theorem says that any even integer $\geq 4$ is a sum of a prime and the product of at most two primes. There are many significant works on Goldbach conjecture (cf. \cite{Vr2} for excellent exposition). For instance, \cite{FL, MM, Ms, MS, Sx, Sj, Wj} are also interesting related works.

Take Table 2.1, (4.5) and (4.6). Denote
\begin{equation}{\msr P}=\mbb Z_+\setminus(\bigcup_{i=1}^8{\msr T}_i). \end{equation}
Then $\{24m+1\mid m\in\msr P\}$ are all the prime of the form $24m+1$ (cf. Table 2.2). We have\psp

{\bf Goldbach Conjecture Modulo 24}\quad $\{20\leq m\in\mbb Z_+\}\subset\{m_1+m_2\mid m_1,m_2\in \msr P\}.$\psp

We have used AI to test the above conjecture  for $m\leq 10^7$ and it holds. But we failed to get a logic proof. The above conjecture implies
that for any integer $m\geq 20$, $24m+2$ is a sum of primes of the forms  $24m_1+1$ and $24m_2+1$. When $m\in \{1, 2, 3, 4, 5, 9, 10, 15, 19\}$,
$24m+2$ is not a sum of primes of the forms  $24m_1+1$ and $24m_2+1$. If $m\in \{6, 7, 8, 11, 12, 13, 14, 16, 17, 18\}$, $24m+2$ is still a sum of primes of the forms  $24m_1+1$ and $24m_2+1$.\psp

\psp

\psp

E-Mail: Y. Cao, caoyichuan@amss.ac.cn; \quad X. Xu, xiaoping@math.ac.cn

 \end{document}